\documentclass[12pt,a4paper,fleqn]{article}

\usepackage{lscape,exscale,amsthm}
\usepackage[intlimits]{amsmath}
\usepackage{rawfonts}
\usepackage{latexsym}
\usepackage[cp850]{inputenc}
\usepackage{epsfig}
\usepackage{bbm}
\usepackage{enumitem}
\usepackage{float}

\usepackage[authoryear]{natbib}
\RequirePackage[OT1]{fontenc}
\RequirePackage{amsthm,amsmath}
\newcommand{\bSigma}{\mathbf{\Sigma}}

\newcommand{\bb}{\mathbf{b}}

\newcommand{\bx}{\mathbf{X}}

\newcommand{\be}{\mathbf{e}}
\newcommand{\by}{\mathbf{Y}}

\newcommand{\ba}{\mathbf{a}}

\newcommand{\bR}{\mathbf{R}}

\newcommand{\bv}{\mathbf{v}}
\newcommand{\bw}{\mathbf{w}}

\newcommand{\bi}{\mathbf{1}}

\newcommand{\bI}{\mathbf{I}}

\newcommand{\bA}{\mathbf{A}}

\newcommand{\bD}{\mathbf{D}}

\newcommand{\bS}{\mathbf{S}}

\newcommand{\sy}{\bar{\mathbf{y}}}

\newcommand{\bu}{\mathbf{u}}
\newcommand{\bU}{\mathbf{U}}
\newcommand{\bTheta}{\mathbf{\Theta}}

\usepackage{tikz}

\usepackage{stackrel}

\numberwithin{equation}{section}
\theoremstyle{plain}
\newtheorem{theorem}{Theorem}[section]

\newtheorem{lemma}{Lemma}[section]
\newtheorem{corollary}{Corollary}[section]

\newtheorem{remark}{Remark}

\begin{document}

\begin{center}
  \vspace*{2cm} \noindent {\bf \large Spectral analysis of the Moore-Penrose inverse of a
  large dimensional
  sample covariance matrix}\\

\vspace{1cm} \noindent {\sc  Taras Bodnar$^{a}$, Holger Dette$^{b}$ and Nestor Parolya$^{c}$}\\
\vspace{1cm}
{\it \footnotesize  $^a$
Department of Mathematics, Stockholm University, SE-10691 Stockhom, Sweden\\
e-mail: taras.bodnar@math.su.se
} \\
{\it \footnotesize  $^b$
Department of Mathematics,  Ruhr University Bochum, D-44870 Bochum, Germany\\
e-mail: holger.dette@ruhr-uni-bochum.de
}\\
{\it \footnotesize  $^c$
Institute of Empirical Economics, Leibniz University Hannover, D-30167 Hannover, Germany \\
e-mail: nestor.parolya@ewifo.uni-hannover.de}

\end{center}
\vspace{1cm}

\begin{abstract}
For a sample of $n$  independent identically distributed $p$-dimensional centered random vectors
with covariance matrix $\bSigma_n$ let $\tilde{\bS}_n$  denote the usual sample covariance
(centered by the mean) and $\bS_n$  the non-centered sample covariance matrix (i.e.
the matrix of second moment estimates), where $p> n$.  In this paper, we provide the limiting spectral distribution and
central limit theorem for linear spectral
statistics of the Moore-Penrose inverse of $\bS_n$ and $\tilde{\bS}_n$.
 We consider the large dimensional asymptotics when the number of variables $p\rightarrow\infty$ and the sample size $n\rightarrow\infty$ such that $p/n\rightarrow c\in (1, +\infty)$. We present a Marchenko-Pastur law for both types of matrices, which shows that the limiting spectral distributions for both sample covariance matrices are the same. On the other hand, we
 demonstrate  that the asymptotic distribution of  linear spectral statistics  of  the Moore-Penrose inverse of $\tilde{\bS}_n$ differs in the
 mean   from that of $\bS_n$.
\end{abstract}

\vspace{0.7cm}

\noindent AMS 2010 subject classifications: 60B20, 60F05, 60F15, 60F17, 62H10\\
\noindent {\it Keywords}: CLT, large-dimensional asymptotics, Moore-Penrose inverse, random matrix theory. \\

\clearpage

\section{Introduction}

Many statistical, financial and genetic problems require estimates of the inverse population covariance matrix which are often constructed by inverting the sample covariance matrix. Nowadays, the modern scientific data sets involve the large number of sample points which is often less than the     dimension  (number of features) and so the sample covariance matrix is not invertible. For example, stock markets include a large number of companies which is often larger than the  number of available time points; or the DNA can contain a fairly large number of genes in comparison to a small number of patients. In such situations, the Moore-Penrose inverse or pseusoinverse of the sample covariance matrix can be used as an estimator for the precision matrix  
  [see, e.g., \cite{srivastava2007},  \cite{kubsriv2008}, \cite{hoyle2011}, \cite{bodguppar2015}].

In order to  better understand the statistical properties of estimators and tests based on  the Moore-Penrose inverse in  high-dimensional settings,
it is of interest to study the asymptotic spectral properties of the Moore-Penrose inverse, for example convergence of its linear spectral statistics (LSS).
 This information is of great interest for high-dimensional statistics because more efficient estimators and tests, which do not suffer from the ``curse of dimensionality'' and do not reduce the number of dimensions, may be constructed and applied in practice.
 Most of the classical multivariate procedures are based on the central limit theorems assuming that the dimension $p$ is fixed and the sample size $n$ increases. However, it has been pointed out by numerous authors that this assumption does not yield precise  distributional approximations for commonly used statistics,
  and that better approximations can be obtained considering scenarios where     the dimension tends to infinity as well [see, e.g., \cite{baisil2004} and references therein]. More precisely, under the high-dimensional asymptotics we understand the case when the sample size $n$ and the dimension $p$ tend to infinity, such that their ratio $p/n$ converges to some positive constant $c$. Under this condition the well-known Marchenko-Pastur equation as well as Marchenko-Pastur law were derived [see, \cite{marpas1967}, \cite{silverstein1995}].

While most authors in random matrix theory  investigate spectral properties of the sample covariance matrix
$ \bS_n=\frac{1}{n} \sum_{i=1}^n \mathbf{y}_i \mathbf{y}_i^\prime$ (here $ \mathbf{y}_1, \ldots ,  \mathbf{y}_n$ denotes a sample of i.i.d.
 $p$-dimensional random vectors  with mean $0$ and variance $\bSigma_n$), \cite{pan2014} studies the differences
occurring if $ \bS_n$ is replaced by its  centered version $ \tilde{\bS}_n=\frac{1}{n} \sum_{i=1}^n ( \mathbf{y}_i-\sy)( \mathbf{y}_i-\sy)^\prime$
(here $\sy$ denotes the mean of $ \mathbf{y}_1, \ldots ,  \mathbf{y}_n$).
Corresponding (asymptotic) spectral properties for the inverse of $ \bS_n$ have been recently derived by
 \cite{zhebaiyao2013} in the case $p < n$, which correspond to the case
$c<1$. The aim of the present  paper is to close a gap in the literature and
focussing on the case $c\in(1, \infty)$. We  investigate the  differences
in the asymptotic spectral properties of  Moore-Penrose inverses of  centered and non-centered sample
covariance matrices.
 In particular  we provide the limiting spectral distribution and the central limit theorem (CLT) for linear spectral statistics of the Moore-Penrose inverse of the sample covariance matrix. \\
 In Section \ref{sec2}  we present the Marchenko-Pastur
  equation together with a Marchenko-Pastur law for the
 Moore-Penrose inverse of the sample covariance matrix. Section \ref{sec3} is divided into two parts: the first one is dedicated to the CLT for the LSS of the pseusoinverse of the non-centered sample covariance matrix while the second part covers the case when the sample covariance matrix is a centered one. 
 While  the limiting spectral distributions for both sample covariance matrices are the same, it is shown that  the asymptotic distribution of  LLS  of  the Moore-Penrose inverse of ${\bS}_n$ and $\tilde{\bS}_n$  differ.
 Finally, some technical details are given in Section \ref{sec4}.

\section{Preliminaries and the Marchenko-Pastur equation }
 \label{sec2}
\def\theequation{2.\arabic{equation}}
\setcounter{equation}{0}

Throughout this paper we use the following notations and assumptions:
\begin{itemize}
\item For a symmetric matrix $\bA$ we denote by $\lambda_1(\bA)\geq\ldots\geq\lambda_p(\bA)$
 its ordered eigenvalues and  by  $F^{\bA}(t)$ the corresponding empirical distribution function (e.d.f.), that is
 $$
 F^{\bA}(t) =   \dfrac{1}{p}\sum\limits_{i=1}^{p}\mathbbm{1}{  \{  \lambda_i(\bA)  \leq t\}}  ,
 $$
 where $\mathbbm{1}{\{\cdot\}}$ is the indicator function.
\item {\bf (A1)} Let $\bx_n$ be a $p\times n$ matrix which consists of independent and identically distributed (i.i.d.) real random variables with zero mean and  unit variance.
\item {\bf (A2)} For the latter matrix $\bx_n=(X_{ij})^{j=1,\ldots,n}_{i=1,\ldots,p}$ we assume additionally that $E(X^{4+\delta}_{11})<\infty$ for some $\delta>0$.

\item By
$$
 \by_n=\bSigma_n^{\frac{1}{2}}\bx_n\,.
$$
 we define a  $p\times n$ observation matrix  with independent columns with mean  $0$ and covariance matrix $\bSigma_n$.
\footnote{{ We could easily include the population mean vector into the model but it will only make the formulas for
weak convergence more complex not the analysis itself.}} It is further assumed that neither $\bSigma_n^{\frac{1}{2}}$ nor $\bx_n$ are observable.

\item  The centered and non-centered sample covariance matrix  are denoted by
 \begin{eqnarray*}\label{samplecov}
 \tilde{\bS}_n&=&\frac{1}{n}(\by_n-\sy\bi^{\prime})(\by_n-\sy\bi^\prime)^\prime=\dfrac{1}{n}\by_n\by_n^{\prime}  - \sy\sy^\prime\, \\
  \label{samplecov2}
 \bS_n&=&\frac{1}{n}\by_n\by_n^\prime=\frac{1}{n}\bSigma_n^{\frac{1}{2}}\bx_n\bx_n^{\prime}\bSigma_n^{\frac{1}{2}}\,,
 \end{eqnarray*}
where $\bi$ denotes the $n$-dimensional vector of ones and $\sy =\frac {1}{n} \sum^n_{i=1} \mathbf{y}_i$. The corresponding e.d.f.'s  are given by  $F^{{{\tilde \bS}_n}}$ and  $F^{\bS_n}$, respectively.

\item The Moore-Penrose inverse of a $p \times p$ matrix $\bA$ is denoted by $\bA^+$ and from its definition must satisfy the following four criteria [see, e.g., \cite{hornjohn1985}]
  \begin{itemize}
  \item[(i)]   $\bA\bA^+\bA=\bA$,
  \item[(ii)]  $\bA^+\bA\bA^+=\bA^+$,
  \item[(iii)] $\bA\bA^+$ is symmetric,
  \item[(iv)]  $\bA^+\bA$ is symmetric.
  \end{itemize}
It is worth pointing out that the generalized inverse considered recently by \cite{bodguppar2015} does not satisfy the conditions \textbf{(iii)} and \textbf{(iv)} presented above. If the matrix $\bA$ has a full column rank then the matrix $\bA^\prime\bA$ is invertible and the Moore-Penrose inverse obeys a simple representation given by
\begin{equation}\label{frmp}
  \bA^+=(\bA^\prime\bA)^{-1}\bA^\prime\,.
\end{equation}

\item
  For  a function $G: \mathbbm{R} \to \mathbbm{R} $ of  bounded variation we introduce
  the Stieltjes transform
$$
  m_G(z)=\int\limits_{-\infty}^{+\infty}\dfrac{1}{\lambda-z}dG(\lambda); ~~~z\in\mathbbm{C}^+ \,.
$$
\end{itemize}

\begin{remark}
{\rm Although   assumption \textbf{(A2)} requires the existence of  moments of order $4+\delta$,   we  suspect  that the results of this paper also hold under the     existence of
  moments of order 4. For a proof one would have to use   truncation   techniques as provided by \cite{baimpan2007} for the matrix $1/n\by^\prime_n\by_n$.  These extremely technical details are omitted for the sake of brevity and transparency.}
\end{remark}
\medskip

 In this paper we are interested in the asymptotic
 properties of the empirical distribution function and linear spectral statistics of the
 eigenvalues of the Moore-Penrose inverse of the matrices $\tilde{\bS}_n$ and $\bS_n$.
 Actually, the limiting spectral distribution of both  matrices coincide because they differ only by a rank
 one perturbation. This is shown in the following Lemma \ref{lem1}.

\begin{lemma} \label{lem1}
Let $\tilde{\bS}^+_n$ and $\bS_n^+$ be the Moore-Penrose inverses  of centered and non-centered sample covariance matrices, respectively, then
\begin{equation} \label{erg1}
  ||F^{\tilde{\bS}^+_n} - F^{\bS^+_n} ||_\infty  ~\leq ~ \frac{2}{p}\,,
\end{equation}
where $|| g ||_\infty$ is the usual supremum norm of a function $g  : \mathbbm{R} \to \mathbbm{R}$.
\end{lemma}

\medskip

{\bf Proof.}  We obtain for the rank one update of  the Moore-Penrose inverse [see \cite{meyer1973}]
\begin{equation}
  \label{eq:moore}
  \tilde{\bS}_n^+ = (\bS_n-\sy\sy^\prime)^+=\bS_n^+-\dfrac{\bS_n^+\sy\sy^\prime(\bS_n^+)^2+(\bS_n^+)^2\sy\sy^\prime(\bS_n^+)}{\sy^\prime(\bS_n^+)^2\sy} +\dfrac{\sy^\prime(\bS_n^+)^3\sy}{(\sy^\prime(\bS_n^+)^2\sy)^2}\bS_n^+\sy\sy^\prime\bS_n^+\,.
\end{equation}
With the notation   $\bu= \bS^+_n\sy$ and $\bv=(\bS^+_n)^2\sy$ the difference of the Moore-Penrose inverses  can
therefore be rewritten as follows
\begin{eqnarray}
\label{diff}
\bS_n^+ -  \tilde{\bS}_n^+&=&\frac{\bu\bv^\prime+\bv\bu^\prime}{\bu^\prime\bu} - \frac{\bu^\prime\bv}{(\bu^\prime\bu)^2}\bu\bu^\prime\nonumber\\
   &=&  \frac{1}{\bu^\prime\bu}\Big(\frac{\bu^\prime\bu}{\bu^\prime\bv}\bv\bv^\prime - \Big( \Big[\frac{\bu^\prime\bv}{\bu^\prime\bu}\Big]^{1/2}\bu - \bv\Big [\frac{\bu^\prime\bu}{\bu^\prime\bv}\Big]^{1/2}  \Big)\Big( \Big[\frac{\bu^\prime\bv}{\bu^\prime\bu}\Big]^{1/2} \bu - \bv \Big[\frac{\bu^\prime\bu}{\bu^\prime\bv}\Big]^{1/2}  \Big)^\prime\Big) \nonumber\\
&=& \frac{1}{\bu^\prime\bv}\bv\bv^\prime - \frac{1}{\bu^\prime\bu}\bw\bw^\prime,
\end{eqnarray}
where
$$
   \bw=\left[\frac{\bu^\prime\bv}{\bu^\prime\bu}\right]^{1/2}\bu - \bv\left [\frac{\bu^\prime\bu}{\bu^\prime\bv}\right]^{1/2}\,.
$$
Thus, the difference $\bS_n^+ -  \tilde{\bS}_n^+$ is a matrix at most of rank $2$ and  the rank inequality in Theorem A.43
of  \cite{baisil2010}   yields the estimate \eqref{erg1}.  \hfill $\Box$

\bigskip

We now study the asymptotic  properties of the e.d.f. of the spectrum of the  Moore-Penrose inverse
of  the  sample covariance matrices $\tilde{\bS}_n$ and $\bS_n$.
 As a consequence of Lemma \ref{lem1} and equality \eqref{frmp}
the asymptotic properties of the e.d.f. of both Moore-Penrose inverses  can be studied
  concentrating  on the matrix
$$
 \bS_n^+= (1/n\by_n\by_n^{\prime})^+ =\left[(1/\sqrt{n}\by_n)^{+}\right]^\prime(1/\sqrt{n}\by_n)^+ = 1/n\by_n(1/n\by_n^\prime\by_n)^{-2}\by_n^{\prime}
 $$
together with the corresponding  e.d.f. $F^{\bS_n^+}$.
Indeed, because
 $m_{F^{1/n\by_n^\prime\by_n}}(z)$ tends almost surely to the solution of the Marchenko-Pastur equation we obtain the following result.

 \begin{theorem}\label{mooreMP} Assume that \textbf{(A1)} holds,  $\frac{p}{n}\rightarrow c\in(1, \infty)$ as $n\rightarrow\infty$ and 
 that $F^{\bSigma_n}$  converges  weakly to a cumulative distribution function  (c.d.f.)
 \ $H$. Then the e.d.f.  $F^{\bS_n^+}$ converges weakly almost surely to some deterministic c.d.f.
  $P$ whose Stieltjes transformation $m_P$  satisfies the following equation
$$
  m_P(z) = -\frac{1}{z}\Bigg (2- c^{-1} + \int_{-\infty}^{+\infty}\dfrac{dH(\tau)}{z\tau c(zm_P(z)+1)-1}\Bigg)\,.
$$
 \end{theorem}

{\bf Proof.}
Let $\bS_n= \bU_n\bD_n\bU_n^\prime $ be the eigenvalue decomposition of the matrix $\bS_n$, then
\begin{eqnarray}
 m_{F^{\bS_n^+}}(z)
 &=& \dfrac{1}{p}\text{tr}\left[\left((\bU_n\bD_n\bU_n^\prime)^+ - z\bI \right)^{-1}\right]
  = \dfrac{1}{p}\text{tr}\left[\left(\bD_n^+ - z\bI \right)^{-1}\right]\nonumber\\
 &=& \dfrac{p-n}{p}\left(-\dfrac{1}{z}\right) + \dfrac{n}{p}\dfrac{1}{n}\sum\limits_{i=1}^n \dfrac{1}{\lambda_i((1/n\by_n^\prime\by_n)^{-1})-z}\nonumber\\
 &=& \left(1-\dfrac{n}{p}\right)\left(-\dfrac{1}{z}\right)+ \dfrac{n}{p}m_{F^{(1/n\by_n^\prime\by_n)^{-1}}}(z)\label{ref1}\,.
\end{eqnarray}
The last two equalities follow from the fact that the spectrum of the matrix $\bS_n^+$ differs from that of $(1/n\by_n^\prime\by_n)^{-1}$ by exactly $p-n$ zero eigenvalues.
For the Stieltjes transform  $m_{F^{(1/n\by_n^\prime\by_n)^{-1}}}(z)$ in this expression we get
\begin{eqnarray}
 m_{F^{(1/n\by_n^\prime\by_n)^{-1}}}(z) &=& \dfrac{1}{n}\sum\limits_{i=1}^n \dfrac{1}{\lambda^{-1}_i(1/n\by_n^\prime\by_n)-z}\nonumber\\
 &=& -\dfrac{1}{nz}\sum\limits_{i=1}^n \dfrac{\lambda_i(1/n\by_n^\prime\by_n)}{\lambda_i(1/n\by_n^\prime\by_n)-\frac{1}{z}}\nonumber\\
 &=& -\dfrac{1}{z} - \dfrac{1}{z^2}\dfrac{1}{n}\sum\limits_{i=1}^n \dfrac{1}{\lambda_i(1/n\by_n^\prime\by_n)-\frac{1}{z}}\nonumber\\
 &=& -\dfrac{1}{z} - \dfrac{1}{z^2}m_{F^{1/n\by_n^\prime\by_n}}\left(\frac{1}{z}\right)\label{ref2}\,.
\end{eqnarray}
Combining the identities (\ref{ref1}) and (\ref{ref2}) provides an equation which relates
the Stieltjes transforms of $F^{\bS_n^+}$ and $ F^{1/n\by^\prime_n\by_n}$, that is
\begin{equation} \label{eq1}
m_{F^{\bS_n^+}}(z) = -\dfrac{1}{z} - \dfrac{n}{p}\dfrac{1}{z^2}m_{F^{1/n\by^\prime_n\by_n}}(1/z)\,.
\end{equation}
 It now follows from   \cite{baisil2010} that as $p/n\rightarrow c>1$ the e.d.f.'s $F^{1/n \by^\prime_n \by_n}$ and $F^{1/n \by_n \by_n^\prime}$ converge weakly almost surely to non-generate distribution functions $\underline{F}$ and $F$ with corresponding Stieltjes transforms satisfying the equation
$$
 m_{\underline{F}}(z)  = -\dfrac{1-c}{z}+cm_{F}(z)\, .
$$
Consequently,  we have from \eqref{eq1} almost surely
  \begin{eqnarray*}
   m_{F^{\bS_n^+}}(z) \longrightarrow  m_P(z)  &:=& -\dfrac{1}{z}- c^{-1}\dfrac{1}{z^2}\left( z(c-1) + cm_{F}(1/z) \right)\\
   &=&  -\dfrac{2-c^{-1}}{z}-\dfrac{m_{F}(1/z)}{z^2}
  \end{eqnarray*}
  as $n\rightarrow\infty$,
  where $m_{F}(1/z)$ is the Stieltjes transform of the limiting distribution $F$ of the e.d.f. of $\bS_n$, which satisfies the  equation (see, \cite{silverstein1995})
  \begin{equation}\label{MPeq}
   m_{F}(1/z) = \int\limits_{-\infty}^{+\infty} \dfrac{dH(\tau)}{\tau(1-c-c\frac{m_{F}(1/z)}{z})-\frac{1}{z}}\,.
  \end{equation}
Thus, $z(c^{-1}-2 - zm_P(z)) = m_{F}(1/z)$ must satisfy the same equation (\ref{MPeq}) which implies
$$
z(c^{-1}-2 - zm_P(z)) = \int\limits_{-\infty}^{+\infty} \dfrac{dH(\tau)}{\tau(1-c-c\frac{z(c^{-1}-2 - zm_P(z))}{z})-\frac{1}{z}}\,.
$$
After some simplification the result follows. \hfill $\Box$

\medskip

\medskip

\begin{corollary}\label{cor1} If $\bSigma_n=\sigma^2\bI_n$ and the assumptions of Theorem \ref{mooreMP} are satisfied,
then the e.d.f. of $\bS_n$ converges weakly almost surely to a deterministic distribution function $P$ with Stieltjes transform
$$
 m_P(z) = -\dfrac{1}{z}\Big(1 + \dfrac{-1/z+(c-1)\sigma^2+\sqrt{(1/z-c\sigma^2+\sigma^2)^2-4/z\sigma^2}}{2\sigma^2c} \Big)\,.
$$
Moreover, the limiting distribution is given by
$$
F= (1-c^{-1})  \delta_0 + \nu(x)dx
$$
where $\delta_a$ denotes the Dirac measure at the point $a \in \mathbbm{R}$ and
\begin{eqnarray*}
\nu(x)=
\left\{
\begin{array}{cr}
c^{-1}\dfrac{\sqrt{(\lambda^{-1}_+-1/x)(1/x-\lambda^{-1}_-)}}{2\pi\sigma^2x}, & x\in[\lambda^{-1}_+,\lambda^{-1}_-]\vspace{5mm}\\
0, & \text{otherwise},
\end{array}
\right.
\end{eqnarray*}
with $\lambda_+=\sigma^2(1+\sqrt{c})^2$ and $\lambda_-=\sigma^2(1-\sqrt{c})^2$.
\end{corollary}

\section{CLT for linear spectral statistics}
 \label{sec3}
\def\theequation{3.\arabic{equation}}
\setcounter{equation}{0}

For a (random) symmetric matrix $\bA$ with spectrum $\lambda_1({\bA}), \ldots, \lambda_p({\bA})$ we consider
the linear spectral statistic
$$
F^{\bA}(g) = \int\limits_{-\infty}^{+\infty}g(x)dF^{\bA }(x)=\frac{1}{p}\sum\limits_{i=1}^p g(\lambda_i(\bA ))
$$
where $g : \mathbbm{R} \to \mathbbm{R}$ is  a given test function. In the next two sections we will investigate
 the asymptotic properties of  $F^{\bS^+_n}$ and $F_n^{\tilde \bS^+_n}$.

\subsection{Linear spectral statistics of $\bS^+_n$}
 \label{sec31}


In the following discussion we consider the random function
$$
  G_n(x) =  p (F^{\bS_n^+}- P^*_n)(x)\,,
$$
where $P^*_n$ is a finite sample proxy of the limiting spectral distribution of $\bS_n^+$, namely $P$. Note that the function $P^*_n$ is constructed simply by substituting $p/n$ for $c$ and $H_n=F^{\bSigma_n}$ for $H$ into the limiting distribution $P$.

 \begin{theorem} \label{cltsimple}

Assume that
  \begin{itemize}
  \item[(i)] \textbf{(A1)} and \textbf{(A2)} hold.
 \item[(ii)] $\dfrac{p}{n}\rightarrow c\in(1, \infty)$ as $n\rightarrow\infty$.
 \item[(iii)] The population covariance matrix $\bSigma_n$ is nonrandom symmetric and positive definite  with a bounded spectral norm and the e.d.f.  $H_n=F^{\bSigma_n}$ converges weakly to a nonrandom distribution function $H$ with a support bounded away from zero.

\item[(iv)]   $g_1,\ldots, g_k$ are functions on $\mathbbm{R}$ analytic on an open region $D$ of the complex plane, which contains the real interval
   $$
     \big [0,~\limsup_{n \to \infty} \lambda_{max}(\bSigma^{-1}_n)/(1-\sqrt{c})^2 \big]
  $$
\item[(v)] Let $\be_i$ denote a $p$-dimensional vector   with the $i$-th element $1$ and others $0$, and define
  \begin{eqnarray*}
\kappa_i (z) & =& \be^\prime_i\bSigma^{1/2}_n(m_{\underline{F} }(1/z)\bSigma_n+\bI)^{-1}\bSigma^{1/2}_n\be_i,\\
    \chi_i(z) &=& \be^\prime_i\bSigma^{1/2}_n(m_{\underline{F} }(1/z)\bSigma_n+\bI)^{-2}\bSigma^{1/2}_n\be_i\,
  \end{eqnarray*}
as the $i$th diagonal elements of the matrices
$$ \bSigma^{1/2}_n(m_{\underline{F} } (1/z)\bSigma_n+\bI)^{-1}\bSigma^{1/2}_n, \quad   \bSigma^{1/2}_n(m_{\underline{F} } (1/z)\bSigma_n+\bI)^{-2}\bSigma^{1/2}_n,$$ respectively. The population covariance matrix $\bSigma_n$
satisfies
\begin{eqnarray*}
 &&  \lim_{n\to \infty}  \frac{1}{n}\sum_{i=1}^n  \kappa_i(z_1) \kappa_i(z_2)  =  h_1(z_1,z_2) \\
&&
  \lim_{n\to \infty}      \frac{1}{n}\sum_{i=1}^n  \kappa_i(z)\chi_i(z) = h_2(z)
  \end{eqnarray*}
  \end{itemize}
  then
$$\left( \int\limits_{-\infty}^{+\infty}g_1(x)dG_n(x),\ldots,  \int\limits_{-\infty}^{+\infty}g_k(x)dG_n(x)\right)'\overset{\mathcal{D}}{\longrightarrow} \left(X_{g_1},\ldots, X_{g_{k}} \right)',   $$
where $\left(X_{g_1},\ldots, X_{g_{k}} \right)'$ is a Gaussian vector with mean
{\small
\begin{eqnarray}
  \label{eq:mean1}
  E(X_{g}) & = & \frac{1}{2\pi i}\oint \frac{g(z)}{z^2}\frac{c\int\limits_{-\infty}^{+\infty}\frac{t^2(m_{\underline{F}}(1/z))^3}
  {(1+tm_{\underline{F}}(1/z))^3}dH(t)}{\Big(1-c\int\limits_{-\infty}^{+\infty}\frac{t^2(m_{\underline{F}}(1/z))^2}{(1+tm_{\underline{F}}(1/z))^2}dH(t)\Big)^2}dz\nonumber\\
&+& \frac{E(X^4_{11})-3}{2\pi i}\oint \frac{g(z)}{z^2}\frac{c(m_{\underline{F}}(1/z))^3h_2(z)}{1-c\int\limits_{-\infty}^{+\infty}\frac{t^2(m_{\underline{F}}(1/z))^2}{(1+tm_{\underline{F}}(1/z))^2}dH(t)}dz\,
\end{eqnarray}
}
and covariance function
{\small
\begin{eqnarray}
  \label{cov1}
{\rm Cov}(X_{g_1}, X_{g_2})& = & -\frac{1}{2\pi^2}\oint\oint\frac{g(z_1)g(z_2)}{z^2_1z^2_2}\frac{m_{\underline{F}}^\prime(1/z_1)m_{\underline{F}}^\prime(1/z_2)}{\left(m_{\underline{F}}(1/z_1)-m_{\underline{F}}(1/z_2)\right)^2}dz_1dz_2\nonumber\\
 &-& \frac{E(X^4_{11})-3}{4\pi^2}\oint \oint\frac{g(z_1)g(z_2)}{z^2_1z^2_2}\left[m_{\underline{F}}(1/z_1)m_{\underline{F}}(1/z_2)h_1(z_1, z_2) \right]''dz_1dz_2\,.
\end{eqnarray}
} The contours in (\ref{eq:mean1}) and  (\ref{cov1}) are both contained in the analytic region for the functions
$g_1, \ldots, g_k$ and both enclose the support of $P^*_n$ for sufficiently large $n$. Moreover, the contours in (\ref{cov1}) are disjoint.
 \end{theorem}

{\bf Proof.}
The proof of this theorem  follows combining arguments from the proof of Lemma 1.1 of \cite{baisil2004} and Theorem 1 of \cite{pan2014}.
      To be precise, we note first that $m_{F^{1/n\by^\prime_n\by_n}}(z)$   converges almost surely   with limit, say $m_{\underline{F}}(z)$. We also observe that the CLT of $p(m_{F^{\bS_n^+}}(z)-m_P(z))$  is the same as of $-\frac{n}{z^2}(m_{F^{1/n\by_n^\prime\by_n}}(1/z)- m_{\underline{F}}(1/z))$ and proceeds in two steps

    \begin{itemize}
    \item[(i)]  Assume first that $E(X^4_{11})=3$.
By Lemma 1.1 of \cite{baisil2004}  we get that the process
$$
  \label{eq:clt2004}
  -\frac{n}{z^2}(m_{F^{1/n\by_n^\prime\by_n}}(1/z)- m_{\underline{F}}(1/z))
$$
defined on an arbitrary positively oriented contour $\mathcal{C}$ which contains the support of $P^*_n$ for sufficiently larger $n$
converges weakly to a Gaussian process with mean
$$
  \label{eq:mean2004}
 - \frac{1}{z^2}\frac{c\int\limits_{-\infty}^{+\infty}\frac{t^2(m_{\underline{F}}(1/z))^3}{(1+tm_{\underline{F}}(1/z))^3}dH(t)}{\Big(1-c\int\limits_{-\infty}^{+\infty}\frac{t^2(m_{\underline{F}}(1/z))^2}{(1+tm_{\underline{F}}(1/z))^2}dH(t)\Big)^2}
$$
and covariance
$$
  \label{eq:covar2004}
  \frac{1}{z^2_1z^2_2}\frac{m_{\underline{F}}^\prime(1/z_1)m_{\underline{F}}^\prime(1/z_2)}{\left(m_{\underline{F}}(1/z_2)-m_{\underline{F}}(1/z_1)  \right)^2}-\frac{1}{(z_1-z_2)^2}\,.
$$
In order to obtain the assertion of Theorem \ref{cltsimple} we use again the argument made at the beginning of the proof and the following identity
\begin{equation}\label{eqcomp}
  \int_{-\infty}^\infty g(x)dG_n(x) = -\frac{1}{2\pi i}\oint g(z)m_{G_n}(z)dz\, ,
\end{equation}
which is valid with probability one for any  analytic function $g$ defined on an open set containing the support of $G_n$ if  $n$ is sufficiently  large.
The complex integral on the RHS of (\ref{eqcomp}) is over certain positively oriented contour enclosing $G$ on which the function $g$ is analytic [see the discussion in \cite{baisil2004} following Lemma 1.1)]. Further, following the proof of Theorem 1.1 in \cite{baisil2004}  we obtain that $ \left( \int\limits g_1(x)dG_n(x),\ldots,  \int g_k(x)dG_n(x)\right)'$ converges weakly to a Gaussian vector $\left(X_{g_1},\ldots, X_{g_{k}} \right)'$.

\item[(ii)]
In the case $E(X^4_{11})\neq3$ we will use a result proved  in \cite{panzhou2008}, more precisely Theorem 1.4 of this reference. Here we find out that in the case $E(X^4_{11})\neq3$ there appears an additional summand in the asymptotic mean and covariance which involve the limiting functions $h_1(z_1,z_2)$ and $h_2(z)$ from assumption \textbf{(v)}, namely for the mean we obtain the additional term
  $$
  \frac{E(X^4_{11})-3}{z^2}\frac{c(m_{\underline{F}}(1/z))^3h_2(z)}{1-c\int\limits_{-\infty}^{+\infty}\frac{t^2(m_{\underline{F}}(1/z))^2}{(1+tm_{\underline{F}}(1/z))^2}dH(t)}
 $$
and for the covariance
$$
 \frac{E(X^4_{11})-3}{z^2_1z^2_2}\left[m_{\underline{F}}(1/z_1)m_{\underline{F}}(1/z_2)h_1(z_1, z_2) \right]''\,.
$$
These new summands arise by studying the limits of products of $E(X^4_{11})-3$ and the diagonal elements of the matrix $\left(\bS_n-1/z\bI\right)^{-1}$
[see, \cite{panzhou2008}, proof of Theorem 1.4]. Imposing the assumption \textbf{(v)} we basically assure that these limits are the same.

 The assertion is now obtained by the same arguments as given in part \textbf{(i)} of this proof.
 \end{itemize}  \hfill $\Box$

 \bigskip

A version of Theorem \ref{cltsimple} has been proved  in \cite{zhebaiyao2013}  for the usual inverse of the sample covariance matrix $\bS^{-1}_n$. It is also not hard to verify the CLT for linear spectral statistics of $\tilde{\bS}^{-1}_n$ in the case $c<1$ with the same limit distribution.
Theorem \ref{cltsimple} shows that in the case $c \geq 1$ there appear additional terms in the asymptotic mean and covariance of the linear spectral statistics corresponding to $\bS^+_n$.

In general, as we can see, the results on the limiting spectral distribution as well as the CLT for linear spectral statistics of $\bS^+_n$ follow more or less from the already known findings which correspond to the matrix $1/n\by_n^\prime\by_n$. In the next section we will show that the  general CLT for LSS of the random matrix  $\tilde{\bS}^+_n$ is different from that of $\bS_n^+$.

\subsection{CLT for linear spectral statistics of $\tilde{\bS}^+_n$.}
 \label{sec32}

The  goal of this section is to show that the CLT for  linear spectral statistics of the Moore-Penrose inverse $\tilde{\bS}_n^+$  differs
 from that of $\bS_n^+$. In order to show why these two CLTs are different   consider again the   identity \eqref{eq:moore}, that is
$$
  \tilde{\bS}_n^+ = (\bS_n-\sy\sy^\prime)^+=\bS_n^+-\dfrac{\bS_n^+\sy\sy^\prime(\bS_n^+)^2+(\bS_n^+)^2\sy\sy^\prime(\bS_n^+)}{\sy^\prime(\bS_n^+)^2\sy} +\dfrac{\sy^\prime(\bS_n^+)^3\sy}{(\sy^\prime(\bS_n^+)^2\sy)^2}\bS_n^+\sy\sy^\prime\bS_n^+\,.
$$
We emphasize here  the difference  to the well known Sherman-Morrison identity for the inverse of a non-singular
matrix $\tilde{\bS}_n$, that is
\begin{equation} \label{fnew}
\tilde{\bS}^{-1}_n=(\bS_n-\sy\sy^\prime)^{-1}=\bS_n^{-1}+\dfrac{\bS_n^{-1}\sy\sy^\prime\bS_n^{-1}}{1-\sy^\prime\bS_n^{-1}\sy}
\end{equation} [see, \cite{shermorr1950}].
Note that it follows from the identity
$$
  1-\sy^\prime\bS^+_n\sy= 1  - \frac{1}{n^2}\bi^\prime\by^\prime(1/\sqrt{n})\by(1/n\by^\prime\by)^{-2}(1/\sqrt{n})\by^\prime\by\bi = 1- \frac{\bi^\prime\bi}{n} =0\,,
$$
  that the right-hand side of the formula \eqref{fnew} is not defined for the Moore-Penrose inverse in the case $p >n$.

 Now consider the resolvent of  $\tilde{\bS}_n^+$, namely
 \begin{eqnarray}
   \label{eq:resolvent}
   \bR(z)&=&\left( \tilde{\bS}_n^+-z\bI\right)^{-1}=\left(\left[ \bS_n-\sy\sy^\prime\right]^+-z\bI\right)^{-1} \nonumber\\
&=& \left(\bS_n^+-\dfrac{\bS_n^+\sy\sy^\prime(\bS_n^+)^2+(\bS_n^+)^2\sy\sy^\prime(\bS_n^+)}{\sy^\prime(\bS_n^+)^2\sy} +\dfrac{\sy^\prime(\bS_n^+)^3\sy}{(\sy^\prime(\bS_n^+)^2\sy)^2}\bS_n^+\sy\sy^\prime\bS_n^+-z\bI \right)^{-1}\nonumber \\
 &=& \left(\bA(z) - \frac{1}{\bu^\prime\bv}\bv\bv^\prime + \frac{1}{\bu^\prime\bu}\bw\bw^\prime   \right)^{-1}\, ,
\nonumber
 \end{eqnarray}
 where we use \eqref{diff}  and the notations
$\bu= \bS^+_n\sy$, $\bv=(\bS^+_n)^2\sy$ and $\bA(z)= \bS_n^+-z\bI$  and
$$
\bw=\left[\frac{\bu^\prime\bv}{\bu^\prime\bu}\right]^{1/2}\bu - \bv\left [\frac{\bu^\prime\bu}{\bu^\prime\bv}\right]^{1/2}\,
$$
to obtain the last identity. A twofold application of the Sherman-Morrison formula  yields the representation
{\small
\begin{eqnarray}\label{resolvent}
  \bR(z)& = & \left(\bA(z) - \frac{1}{\bu^\prime\bv}\bv\bv^\prime\right)^{-1} -\frac{\left(\bA(z) - \frac{1}{\bu^\prime\bv}\bv\bv^\prime\right)^{-1} \frac{\bw\bw^\prime}{\bu^\prime\bu}\left(\bA(z) - \frac{1}{\bu^\prime\bv}\bv\bv^\prime\right)^{-1}}{1+ \frac{1}{\bu^\prime\bu}\bw^{\prime}\left(\bA(z) - \frac{1}{\bu^\prime\bv}\bv\bv^\prime\right)^{-1}\bw}\nonumber\\
&=& \bA^{-1}(z) + \frac{\frac{1}{\bu^\prime\bv}\bA^{-1}(z)\bv\bv^\prime\bA^{-1}(z)}{1-\frac{1}{\bu^\prime\bv}\bv^\prime\bA^{-1}(z)\bv} -\frac{\left(\bA(z) - \frac{1}{\bu^\prime\bv}\bv\bv^\prime\right)^{-1} \frac{\bw\bw^\prime}{\bu^\prime\bu}\left(\bA(z) - \frac{1}{\bu^\prime\bv}\bv\bv^\prime\right)^{-1}}{1+ \frac{1}{\bu^\prime\bu}\bw^{\prime}\left(\bA(z) - \frac{1}{\bu^\prime\bv}\bv\bv^\prime\right)^{-1}\bw}.
\end{eqnarray}
}
Taking the trace of both sides of (\ref{resolvent}) we obtain the identity
\begin{eqnarray}\label{differ}
 p m_{F^{\tilde{\bS}_n^+}}(z) &=& \text{tr}(\bR(z))  =  p m_{F^{\bS_n^+}}(z)\nonumber\\
 & +&{ \frac{\bv^\prime\bA^{-2}(z)\bv}{\bu^\prime\bv-\bv^\prime\bA^{-1}(z)\bv}-\frac{\bw^\prime\left(\bA^{-1}(z)+\frac{\bA^{-1}(z)
 \bv\bv^\prime\bA^{-1}(z)}{\bu^\prime\bv-\bv^\prime\bA^{-1}(z)\bv} \right)^2\bw}{\bu^\prime\bu+\bw^\prime\bA^{-1}(z)\bw+\frac{(\bw^\prime\bA^{-1}(z)\bv)^2}{\bu^\prime\bv-\bv^\prime\bA^{-1}(z)\bv}} }\, ,
\end{eqnarray}
which  indicates that the CLTs for linear spectral statistics of Moore-Penrose sample covariance matrices of $\bS^+_n$ and $\tilde{\bS}^+_n$
might  differ. In fact, the following result shows that the last two terms on the right-hand side of (\ref{differ}) are asymptotically not negligible.

\begin{theorem}\label{mainteo} Let $ \tilde{G}_n(x) =  p (F^{\tilde{\bS}_n^+}(g)- P^*_n(g))$ and suppose that the assumptions of Theorem \ref{cltsimple} are satisfied,
  then
$$\Big( \int\limits_{-\infty}^{+\infty}g_1(x)d\tilde{G}_n(x),\ldots,  \int\limits_{-\infty}^{+\infty}g_k(x)d\tilde{G}_n(x)\Big)'\overset{\mathcal{D}}{\longrightarrow} \left(X_{g_1},\ldots, X_{g_{k}} \right)',   $$
where $\left(X_{g_1},\ldots, X_{g_{k}} \right)'$ is a Gaussian vector with mean
{\footnotesize
\begin{eqnarray}
  \label{eq:mean}
  E(X_{g}) & = & -\frac{1}{2\pi i}\oint \frac{g(z)}{z^2}\frac{c\int\limits_{-\infty}^{+\infty}\frac{t^2(m_{\underline{F}}(1/z))^3}{(1+tm_{\underline{F}}(1/z))^3}dH(t)}{\big(1-c\int\limits_{-\infty}^{+\infty}\frac{t^2(m_{\underline{F}}(1/z))^2}{(1+tm_{\underline{F}}(1/z))^2}dH(t)\big)^2}dz- \frac{E(X^4_{11})-3}{2\pi i}\oint \frac{g(z)c(m_{\underline{F}}(1/z))^3h_2(z)}{1-c\int\limits_{-\infty}^{+\infty}\frac{t^2(m_{\underline{F}}(1/z))^2}{(1+tm_{\underline{F}}(1/z))^2}dH(t)}dz\nonumber\\
 & - & \frac{1}{2\pi i}\oint \frac{g(z)}{z^2}\frac{m^\prime_{\underline{F}}(1/z)}{m_{\underline{F}}(1/z)}dz
\end{eqnarray}
}
and covariance function
{\small
\begin{eqnarray}
  \label{eq:variance}
\text{Cov}(X_{g_1}, X_{g_2})& = & -\frac{1}{2\pi^2}\oint\oint\frac{g(z_1)g(z_2)}{z^2_1z^2_2}\frac{m_{\underline{F}}^\prime(1/z_1)m_{\underline{F}}^\prime(1/z_2)}{\left(m_{\underline{F}}(1/z_2)\right)^2}dz_1dz_2\nonumber\\
 &-& \frac{E(X^4_{11})-3}{4\pi^2}\oint \oint\frac{g(z_1)g(z_2)}{z^2_1z^2_2}\left[m_{\underline{F}}(1/z_1)m_{\underline{F}}(1/z_2)h_1(z_1, z_2) \right]''dz_1dz_2\,.
\end{eqnarray}}
The contours in (\ref{eq:mean}) and  (\ref{eq:variance}) are both contained in the analytic region for the functions
$g_1, \ldots, g_k$ and both enclose the support of $P_n^*$ for sufficiently large $n$. Moreover, the contours in (\ref{eq:variance}) are disjoint.
\end{theorem}

Before we provide the proof of this result we emphasize the existence of  the extra summand in the asymptotic mean. It has a very simple structure and can be calculated without much   effort in practice.
Indeed, the last integral in \eqref{eq:mean} can be rewritten using integration by parts in the following way (see, Section 4.1 for detailed derivation)
\begin{equation}\label{intcalc}
  - \frac{1}{2\pi i}\oint \frac{g(z)}{z^2}\frac{m^\prime_{\underline{F}}(1/z)}{m_{\underline{F}}(1/z)}dz =  -\frac{1}{\pi }\int_a^b g^\prime(x)\arg[m_{\underline{F}}(1/x)]dx\,,
\end{equation}
where $m_{\underline{F}}(1/x)\equiv \lim\limits_{z\rightarrow x}m_{\underline{F}}(1/z)$ for $x\in\mathbbm{R}$ and the interval $(a, b)$ contains the support of $P$.

On the other hand, the asymptotic variances of the linear spectral statistics for centered and non-centered sample covariance matrices coincide. For a discussion of     assumption $\mathbf{(v)}$ we refer to \cite{pan2014}, Remark 2, and other references therein.

\medskip

\medskip

{\bf Proof of Theorem \ref{mainteo}.}
The proof of Theorem \ref{mainteo} is  based on the Stieltjes transform method and consists of a combination of   arguments similar to those   given by \cite{baisil2004}  and \cite{pan2014}.
First, by the analyticity of  the functions $g_1,\ldots,g_k$ and \eqref{eqcomp} it is sufficient to consider the Stieltjes transforms of the spectral e.d.f. of sample covariance matrices.
Furthermore, recall from \eqref{differ} that the Stieltjes transform of
the e.d.f. of $\tilde{\bS}^+_n$ can be decomposed  as the sum of the Stieltjes transform of
the e.d.f. of $\bS^+_n$ and
the additional term
\begin{equation} \label{xin}
\xi_n(z)
={ \frac{\bv^\prime\bA^{-2}(z)\bv}{\bu^\prime\bv-\bv^\prime\bA^{-1}(z)\bv} -
\frac{\bw^\prime\Big(\bA^{-1}(z)+\frac{\bA^{-1}(z)\bv\bv^\prime\bA^{-1}(z)}{\bu^\prime\bv-\bv^\prime\bA^{-1}(z)\bv} \Big)^2
\bw}{\bu^\prime\bu+\bw^\prime\bA^{-1}(z)\bw+\frac{(\bw^\prime\bA^{-1}(z)\bv)^2}{\bu^\prime\bv-\bv^\prime\bA^{-1}(z)\bv}} }
\end{equation}
 involving sample mean $\sy$ and $\bS^+_n$.
Thus, it is sufficient to show  that this random variable  converges almost surely on
$\mathbbm{C}^+=\{z\in\mathbbm{C}: \Im z>0\}$ to  a nonrandom quantity  as $p/n\rightarrow c>1$ and to determine its limit.
As a result, Theorem \ref{mainteo}  follows from
Slutsky's theorem, the continuous mapping theorem and the results in \cite{baisil2004}  and \cite{panzhou2008}.

It is shown in Section \ref{sec4.2} that the function $\xi_n$ in \eqref{xin} can be represented as
\begin{equation}
  \label{xi_simple}
\xi_n(z)  = -\frac{1}{z}- \frac{\sy^\prime\sy+2z\theta_n(z)+z^2\theta_n^\prime(z)}{1+z\sy^\prime\sy+z^2\theta_n(z)}  \,,
\end{equation}
where the functions $\theta_n $ is given by
\begin{eqnarray} \label{thetrep}
 && \theta_n(z) =  -\frac{1}{z}\sy^\prime\sy +\frac{1}{z}\frac{1}{n}\bi_n^\prime ((1/n\by^\prime_n\by_n)^{-1}-z\bI)^{-1}\bi_n. 
 \end{eqnarray}
As a consequence, the asymptotic properties of $\xi_n$  can be obtained analyzing the quantity
\begin{eqnarray*}
\eta_n(z)& =& \frac{1}{n}\bi_n^\prime ((1/n\by^\prime_n\by_n)^{-1}-z\bI)^{-1}\bi_n\\
&=&  -\frac{1}{z} \text{tr}\left[(1/n\by^\prime_n\by_n) (1/n\by^\prime_n\by_n-1/z\bI)^{-1}1/n\bi_n\bi_n^\prime\right]\\
&=&  -\frac{1}{z} -\frac{1}{z^2}\text{tr}\left[(1/n\by^\prime_n\by_n-1/z\bI)^{-1}\bTheta_n\right]\,,
\end{eqnarray*}
where we use the notation  $\bTheta_n=1/n\bi_n\bi_n^\prime$. It now follows from Theorem 1
 in \cite{rubmes2011}  that
$$
  \left|\text{tr}\left[(1/n\by^\prime_n\by_n-1/z\bI)^{-1}\bTheta_n\right]-x_n(1/z)\right|\longrightarrow0~~\text{a.s.}\, ,
  $$
 where $x_n(1/z)$ is a unique solution in $\mathbbm{C}^+$ of the equation
$$
\dfrac{1+1/zx_n(1/z)}{x_n(1/z)}=\dfrac{c}{p}\text{tr}\left(x_n(1/z)\bI+\bSigma^{-1}_n\right)^{-1}\, .
$$
Note that
  $\text{tr}\left(\mathbf{\Theta}_n\right)=1$ and that
Theorem 1 in \cite{rubmes2011} is originally proven assuming the existence of moments of order $8+\delta$. However, it is shown in    \cite{bodguppar2015} that  only  the existence of moments of order $4+\delta$  is required for this statement.

In order to see how $x_n(1/z)$ relates to $m_{\underline{F}}(1/z)$ we note that due to assumption $\mathbf{(iii)}$   $H_n \stackrel{\mathcal{D}}{\rightarrow} H$ as $n\rightarrow\infty$ and, thus, $x_n(1/z)\rightarrow x(1/z)$. This implies
$$
  \dfrac{1}{x_n(1/z)}+1/z=\dfrac{c}{p}\text{tr}\left(x_n(1/z)\bI+\bSigma^{-1}_n\right)^{-1} =  c \int_{-\infty}^\infty \frac{\tau dH_n(\tau)}{x_n(1/z)\tau+1}\longrightarrow c \int_{-\infty}^\infty \frac{\tau dH(\tau)}{x(1/z)\tau+1},
$$
which leads to
$$
  x(1/z) = -\Big(1/z- c \int_{-\infty}^\infty \frac{\tau}{x(1/z)\tau+1}dH(\tau)  \Big)^{-1}\,.
$$
The last equation is the well-known  Marchenko-Pastur equation for the  Stieltjes transformation ${m}_{\underline F}(1/z)$ of the limiting distribution $\underline{F}$. Because the solution of this equation is  unique  we obtain
$$
  x(1/z) =  {m}_{\underline{F}}(1/z)\,.
$$
As a result we get the following asymptotics for $1+z\sy^\prime\sy+z^2\theta_n(z)$ as $n\rightarrow\infty$
$$
  1+z\sy^\prime\sy+z^2\theta_n(z) \longrightarrow -\frac{m_{\underline{F}}(1/z)}{z}~~\text{a.s.} \,,
$$
which ensures that
\begin{equation}\label{as_xi}
  \xi_n(z) \longrightarrow -\frac{1}{z}-\frac{\left(-\frac{m_{\underline{F}}(1/z)}{z}\right)^\prime}{-\frac{m_{\underline{F}}(1/z)}{z}}=-\frac{1}{z}-z\frac{-\frac{m^\prime_{\underline{F}}(1/z)}{z^3}-\frac{m_{\underline{F}}(1/z)}{z^2}}{m_{\underline{F}}(1/z)}=\frac{1}{z^2}\frac{m^\prime_{\underline{F}}(1/z)}{m_{\underline{F}}(1/z)}~~\text{a.s.}
\end{equation}
for $p/n\rightarrow c \in(1,+\infty)$ as $n\rightarrow\infty$. The assertion of Theorem \ref{mainteo} now  follows
taking into account the argument made at the beginning of the proof and \eqref{as_xi}. \hfill $\Box$

\section{Appendix}
\label{sec4}
\def\theequation{4.\arabic{equation}}
\setcounter{equation}{0}

\subsection{Derivation of \eqref{intcalc}}\label{sec4.1}

We select the contour to be a rectangle with sides parallel to the axes. It intersects the real axis at the points $a\neq 0$ and $b$ such that the interval $(a, b)$ contains the support of $P$. The horizontal sides are taken as a distance $y_0>0$ from the real axis. More precisely,   the contour $\mathcal{C}$ is given by
\begin{equation}
  \mathcal{C} = \big\{a+iy: |y|\leq y_0 \big\} \cup \big\{x+iy_0: x\in [a, b]\big\} \cup \big\{b+iy: |y|\leq y_0 \big\} \cup \big\{x-iy_0: x\in [a, b]\big\}
\end{equation}
so that $(a, b)$ contains $\big [0,~\limsup_{n \to \infty} \lambda_{max}(\bSigma^{-1}_n)/(1-\sqrt{c})^2 \big]$ and is enclosed in the analytic region of function $g$. We calculate the four parts of the contour integral and then let $y_0$ tend to zero. First, we note that
\begin{equation*}
\frac{d}{dz}\log(m_{\underline{F}}(1/z)) = -\frac{1}{z^2}\frac{m^\prime_{\underline{F}}(1/z)}{m_{\underline{F}}(1/z)}\,.
\end{equation*}
Then using integration by parts the last integral in \eqref{eq:mean} becomes
\begin{eqnarray}\label{integral}
   -\frac{1}{2\pi i}\oint \frac{g(z)}{z^2}\frac{m^\prime_{\underline{F}}(1/z)}{m_{\underline{F}}(1/z)}dz& =& \frac{1}{2\pi i}\oint g(z)d(\log m_{\underline{F}}(1/z))\nonumber\\
&=& -\frac{1}{2\pi i}\oint g^\prime(z)\log m_{\underline{F}}(1/z)dz\nonumber\\
&=&-\frac{1}{2\pi i}\oint g^\prime(z)\left(\log|m_{\underline{F}}(1/z)|+i\arg(m_{\underline{F}}(1/z))\right)dz  \,,
\end{eqnarray}
where any branch of the logarithm may be taken. Naturally extending the Stieltjes transform on the negative imaginary axis and using $z=x+iy$ we get
\begin{equation}\label{stbound}
 \big|m_{\underline{F}}(1/z)\big| = \left|\int_{-\infty}^\infty \frac{d\underline{F}(\lambda)}{\lambda-1/z}\right|\leq \int_{-\infty}^\infty \frac{d\underline{F}(\lambda)}{|\lambda-1/z|} =
 \int_{-\infty}^\infty \frac{d\underline{F}(\lambda)}{\sqrt{(\lambda-\frac{x}{|z|^2})^2+\frac{y^2}{|z|^4}}}
 \leq \frac{x^2+y^2}{|y|} \,.
\end{equation}
Next we note that any portion of the integral \eqref{integral} which involves the vertical side can be neglected. Indeed, using \eqref{stbound}, the fact that $|g^\prime(z)|\leq K$ and (5.1) in \cite{baisil2004} for the left vertical side we have
\begin{eqnarray}
&& \left| \frac{1}{2\pi i}\int_{-y_0}^{y_0}g^\prime(a+iy)\log m_{\underline{F}}(1/(a+iy))dy\right|\nonumber\\
&\leq& \frac{K}{2\pi} \int_{-y_0}^{y_0}\left(\log \big|m_{\underline{F}}(1/(a+iy))\big| +\big|\arg(m_{\underline{F}}(1/(a+iy)))\big|\right)dy\nonumber\\
&\leq& \frac{K}{2\pi} \int_{-y_0}^{y_0}\left(\log\frac{a^2+y^2}{|y|} + \pi\right)dy\nonumber\\
&=& \frac{K}{\pi} \int_{0}^{y_0}\log\frac{a^2+y^2}{y}dy +  Ky_0\nonumber\\
&=& \frac{K}{\pi} \left(y_0\log\frac{y^2_0+a^2}{y_0}-y_0+2a\arctan(y_0/a)   \right) +  Ky_0 \label{ibound}\, , 
\end{eqnarray}
 which converges to zero as $y_0\rightarrow0$. A similar argument can be used for the right vertical side. Consequently, only the remaining terms  
\begin{eqnarray*}
 &-& \frac{1}{2\pi}  \int_a^b \Im\left[g^\prime(x+iy_0)\right]\log|m_{\underline{F}}(1/(x+iy_0))|dx - \frac{1}{2\pi}\int_a^b \Re\left[g^\prime(x+iy_0)\right]\arg\bigl{[}m_{\underline{F}}(1/(x+iy_0)) \bigr{]}dx  \\
&-& \frac{1}{2\pi}  \int_a^b \Im\left[g^\prime(x-iy_0)\right]\log|m_{\underline{F}}(1/(x-iy_0))|dx - \frac{1}{2\pi}\int_a^b \Re\left[g^\prime(x-iy_0)\right]\arg\bigl{[}m_{\underline{F}}(1/(x-iy_0)) \bigr{]}dx\, 
\end{eqnarray*}
have to be considered in the limit of the   integral   \eqref{integral}.
Similarly, using the fact that  
\begin{equation}\label{gbound}
  \sup\limits_{x\in[a, b]}|\Im h(x+iy)|
  \leq K|y|\,.
\end{equation}
 for any real-valued analytic function $h$ on the bounded interval $[a, b]$ [see equation (5.6) in \cite{baisil2004}],
 we obtain that the first and the third integrals  are bounded in absolute value by $O(y_0\log y_0^{-1})$ and, thus, can be neglected. As a result the dominated convergence theorem leads to
\begin{equation}
  -\frac{1}{2\pi i}\oint g^\prime(z)\left(\log|m_{\underline{F}}(1/z)|+i\arg(m_{\underline{F}}(1/z))\right)dz \underset{y_0\rightarrow0}{\longrightarrow} - \frac{1}{\pi}\int_a^b g^\prime(x)\arg\bigl{[}m_{\underline{F}}(1/x) \bigr{]}dx\,,
\end{equation}
which proves \eqref{intcalc}.

\subsection{Proof of the representations (\ref{xi_simple}) and (\ref{thetrep})}\label{sec4.2}
For a proof of  (\ref{xi_simple}) we introduce the notations $$\ba=\frac{1}{(\bu^\prime\bv)^{1/2}}\bv~,~~\bb=\frac{(\bu^\prime\bv)^{1/2}}{\bu^\prime\bu}\bu~.
$$ Then $\frac{1}{(\bu^\prime\bu)^{1/2}}\bw=\bb-\ba$, and we obtain for (\ref{xin})
{\footnotesize
\begin{eqnarray*}
  \xi_n(z)& =& \frac{\ba^\prime\bA^{-2}(z)\ba}{1-\ba^\prime\bA^{-1}(z)\ba}\nonumber\\
&-&\frac{(\bb-\ba)^\prime\left(\bA^{-2}(z) + \dfrac{\bA^{-2}(z)\ba\ba^\prime\bA^{-1}(z)+\bA^{-1}(z)\ba\ba^\prime\bA^{-2}(z)}{1-\ba^\prime\bA^{-1}(z)\ba}+\dfrac{\ba^\prime\bA^{-2}(z)\ba (\bA^{-1}(z)\ba\ba^\prime\bA^{-1}(z))}{(1-\ba^\prime\bA^{-1}(z)\ba)^2}  \right) (\bb-\ba)}{1+(\bb-\ba)^\prime\bA^{-1}(z)(\bb-\ba)+\dfrac{((\bb-\ba)^\prime\bA^{-1}(z)\ba)^2}{1-\ba^\prime\bA^{-1}(z)\ba}}\,.
\end{eqnarray*}
}
A tedious but straightforward calculation now gives
{\footnotesize
\begin{eqnarray}\label{xi}
 \xi_n(z) &=& \frac{\ba^\prime\bA^{-2}(z)\ba+\ba^\prime\bA^{-2}(z)\ba(\bb-\ba)^\prime\bA^{-1}(z)(\bb-\ba)+\frac{\ba^\prime\bA^{-2}(z)\ba((\bb-\ba)^\prime\bA^{-1}(z)\ba)^2}{1-\ba^\prime\bA^{-1}(z)\ba}}{1-\ba^\prime\bA^{-1}(z)\ba+(1-\ba^\prime\bA^{-1}(z)\ba)(\bb-\ba)^\prime\bA^{-1}(z)(\bb-\ba)+((\bb-\ba)^\prime\bA^{-1}(z)\ba)^2}   \nonumber\\[0.5cm]
&-&\frac{(1-\ba^\prime\bA^{-1}(z)\ba)(\bb-\ba)^\prime\bA^{-2}(z)(\bb-\ba)}{1-\ba^\prime\bA^{-1}(z)\ba+(1-\ba^\prime\bA^{-1}(z)\ba)(\bb-\ba)^\prime\bA^{-1}(z)(\bb-\ba)+((\bb-\ba)^\prime\bA^{-1}(z)\ba)^2}\nonumber\\[0.5cm]
&-&\frac{2(\bb-\ba)^\prime\bA^{-2}(z)\ba(\bb-\ba)^\prime\bA^{-1}(z)\ba+\frac{\ba^\prime\bA^{-2}(z)\ba((\bb-\ba)^\prime\bA^{-1}(z)\ba)^2}{1-\ba^\prime\bA^{-1}(z)\ba}}{1-\ba^\prime\bA^{-1}(z)\ba+(1-\ba^\prime\bA^{-1}(z)\ba)(\bb-\ba)^\prime\bA^{-1}(z)(\bb-\ba)+((\bb-\ba)^\prime\bA^{-1}(z)\ba)^2}\nonumber\\[0.5cm]
&=& \frac{\ba^\prime\bA^{-2}(z)\ba[\bb^\prime\bA^{-1}(z)\bb]-\bb^\prime\bA^{-2}(z)\bb+2\ba^\prime\bA^{-2}(z)\bb
}{1+\bb^\prime\bA^{-1}(z)\bb-2\ba^\prime\bA^{-1}(z)\bb-\ba^\prime\bA^{-1}(z)\ba[\bb^\prime\bA^{-1}(z)\bb]+(\ba^\prime\bA^{-1}(z)\bb)^2}\nonumber\\[0.5cm]
&+& \frac{\ba^\prime\bA^{-1}(z)\ba[\bb^\prime\bA^{-2}(z)\bb]-2\ba^\prime\bA^{-2}(z)\bb[\ba^\prime\bA^{-1}(z)\bb]}{1+\bb^\prime\bA^{-1}(z)\bb-2\ba^\prime\bA^{-1}(z)\bb-\ba^\prime\bA^{-1}(z)\ba[\bb^\prime\bA^{-1}(z)\bb]+(\ba^\prime\bA^{-1}(z)\bb)^2}\nonumber\\[0.5cm]
&=& \dfrac{\ba^\prime\bA^{-2}(z)\ba[\bb^\prime\bA^{-1}(z)\bb] - \bb^\prime\bA^{-2}(z)\bb [1-\ba^\prime\bA^{-1}(z)\ba]+2\ba^\prime\bA^{-2}(z)\bb[1-\ba^\prime\bA^{-1}(z)\bb] }{[\ba^\prime\bA^{-1}(z)\bb-1]^2+\bb^\prime\bA^{-1}(z)\bb[1-\ba^\prime\bA^{-1}(z)\ba]}\,.
\end{eqnarray}}

Now note that $\xi_n(z)$ is a non-linear function of the quantities $\ba^\prime\bA^{-1}(z)\ba$, $\ba^\prime\bA^{-2}(z)\ba$, $\bb^\prime\bA^{-1}(z)\bb$,
$\bb^\prime\bA^{-2}(z)\bb$, $\ba^\prime\bA^{-1}(z)\bb$ and $\ba^\prime\bA^{-2}(z)\bb$, which can easily  be expressed in term of $\sy$ and $\bS^+_n$. For example, we have
\begin{eqnarray*}
  \ba^\prime\bA^{-1}(z)\ba& =&\frac{1}{\sy^\prime(\bS^+_n)^3\sy}\sy^\prime(\bS_n^+)^2(\bS_n^+-z\bI)^{-1}(\bS_n^+)^2\sy\\
 &=& \frac{1}{\sy^\prime(\bS^+_n)^3\sy}\left(  \sy^\prime\bS_n^+(\bS_n^+-z\bI+z\bI)(\bS_n^+-z\bI)^{-1}(\bS_n^+)^2\sy \right)\\
&=& \frac{1}{\sy^\prime(\bS^+_n)^3\sy}\left(\sy^\prime(\bS^+_n)^3\sy+ z \sy^\prime(\bS_n^+-z\bI+z\bI)(\bS_n^+-z\bI)^{-1}(\bS_n^+)^2\sy  \right)\\
&=& \frac{1}{\sy^\prime(\bS^+_n)^3\sy}\left(\sy^\prime(\bS^+_n)^3\sy + z\sy^\prime(\bS^+_n)^2\sy+ z^2\sy^\prime(\bS_n^+-z\bI)^{-1}(\bS_n^+-z\bI+z\bI)\bS_n^+\sy  \right)\\
&=& \frac{1}{\sy^\prime(\bS^+_n)^3\sy}\left(\sy^\prime(\bS^+_n)^3\sy + z\sy^\prime(\bS^+_n)^2\sy+ z^2\sy^\prime\bS^+_n\sy+
z^3\sy^\prime(\bS_n^+-z\bI)^{-1}(\bS_n^+-z\bI+z\bI)\sy  \right)\\
&=& \frac{1}{\sy^\prime(\bS^+_n)^3\sy}\left(\sy^\prime(\bS^+_n)^3\sy + z\sy^\prime(\bS^+_n)^2\sy+ z^2+z^3\sy^\prime\sy+
z^4\sy^\prime(\bS_n^+-z\bI)^{-1}\sy  \right),
\end{eqnarray*}
where
the last equality follows from the fact that $\sy^\prime\bS^+_n\sy=1$. In a similar way, namely adding and subtracting $z\bI$ from $\bS_n^+$ in a sequel we  obtain representations for the remaining quantities of interest
\begin{eqnarray*}
  \ba^\prime\bA^{-2}(z)\ba & =&  \frac{1}{\sy^\prime(\bS^+_n)^3\sy} \left(\sy^\prime(\bS^+_n)^2\sy+ 2z +  3z^2\sy^\prime\sy + 4z^3\sy^\prime(\bS_n^+-z\bI)^{-1}\sy + z^4\sy^\prime(\bS_n^+-z\bI)^{-2}\sy \right)\\
 \bb^\prime\bA^{-1}(z)\bb &=& \frac{\sy^\prime(\bS^+_n)^3\sy}{(\sy^\prime(\bS^+_n)^2\sy)^2}\left( 1+ z\sy^\prime\sy+ z^2\sy^\prime(\bS_n^+-z\bI)^{-1}\sy\right)\\
 \bb^\prime\bA^{-2}(z)\bb &=& \frac{\sy^\prime(\bS^+_n)^3\sy}{(\sy^\prime(\bS^+_n)^2\sy)^2}\left( \sy^\prime\sy+2z\sy^\prime(\bS_n^+-z\bI)^{-1}\sy+z^2\sy^\prime(\bS_n^+-z\bI)^{-2}\sy\right)\\
\ba^\prime\bA^{-1}(z)\bb &=&  \frac{1}{\sy^\prime(\bS^+_n)^2\sy} \left(\sy^\prime(\bS^+_n)^2\sy+ z+z^2\sy^\prime\sy+z^3\sy^\prime(\bS_n^+-z\bI)^{-1}\sy \right) \\
\ba^\prime\bA^{-2}(z)\bb &=&\frac{1}{\sy^\prime(\bS^+_n)^2\sy} \left(1 +2z\sy^\prime\sy+3z^2\sy^\prime(\bS_n^+-z\bI)^{-1}\sy+z^3\sy^\prime(\bS_n^+-z\bI)^{-2}\sy \right) \,.
\end{eqnarray*}
In the next step we substitute these results in \eqref{xi}. More precisely  introducing
the notations
$$\theta_n(z)=\sy^\prime(\bS_n^+-z\bI)^{-1}\sy, ~\alpha_n=1/\sy^\prime(\bS_n^+)^2\sy,$$
we   have $\theta_n^\prime(z)=\frac{\partial}{\partial z}\theta_n(z)=\sy^\prime(\bS_n^+-z\bI)^{-2}\sy$ and
obtain  for $\xi_n(z)$ the following representation
{\small
\begin{eqnarray*}
 \xi_n(z) &=&\frac{\alpha_n^2(\alpha_n^{-1}+2z+3z^2\sy^\prime\sy+4z^3\theta_n(z)+z^4\theta_n^\prime(z))(1+z\sy^\prime\sy+z^2\theta_n(z))}
 {\alpha_n^2z^2(1+z\sy^\prime\sy+z^2\theta_n(z))^2-\alpha_n^2z(\alpha_n^{-1}+z+z^2\sy^\prime\sy+z^3\theta_n(z))(1+z\sy^\prime\sy+z^2\theta_n(z))}\nonumber\\[0.5cm]
&+&\frac{\alpha_n^2z(\sy^\prime\sy+2z\theta_n(z)+z^2\theta_n^\prime(z))(\alpha_n^{-1}+z+z^2\sy^\prime\sy+z^3\theta_n(z))}
{\alpha_n^2z^2(1+z\sy^\prime\sy+z^2\theta_n(z))^2-\alpha_n^2z(\alpha_n^{-1}+z+z^2\sy^\prime\sy+z^3\theta_n(z))(1+z\sy^\prime\sy+z^2\theta_n(z))}\nonumber\\[0.5cm]
&-&\frac{2\alpha_n^2z(1+2z\sy^\prime\sy+3z^2\theta_n(z)+z^3\theta_n^\prime(z))(1+z\sy^\prime\sy+z^2\theta_n(z))}
{\alpha_n^2z^2(1+z\sy^\prime\sy+z^2\theta_n(z))^2-\alpha_n^2z(\alpha_n^{-1}+z+z^2\sy^\prime\sy+z^3\theta_n(z))(1+z\sy^\prime\sy+z^2\theta_n(z))}\nonumber\\[0.5cm]
&=& \frac{\alpha^2_n\left( (1+z\sy^\prime\sy+z^2\theta_n(z))(\alpha_n^{-1}-z^2(\sy^\prime\sy+2z\theta_n(z)+z^2\theta_n^\prime(z)))\right)}{-\alpha_nz(1+z\sy^\prime\sy+z^2\theta_n(z))}\nonumber\\[0.5cm]
&+&  \frac{\alpha^2_nz(\sy^\prime\sy+2z\theta_n(z)+z^2\theta_n^\prime(z))(\alpha_n^{-1}+z(1+z\sy^\prime\sy+z^2\theta_n(z)))}{-\alpha_nz(1+z\sy^\prime\sy+z^2\theta_n(z))}.
\end{eqnarray*}
}
In order to simplify the following calculations we introduce the quantities
\begin{eqnarray*}
\psi_n(z) = 1+z\sy^\prime\sy+z^2\theta_n(z) ; \qquad
  \psi^\prime_n(z)  =  \sy^\prime\sy + 2z\theta_n(z)+z^2\theta_n^\prime(z)\, ,
\end{eqnarray*}
 which lead to
\begin{eqnarray*}
  \xi_n(z)&=& \frac{\alpha_n^2\psi_n(z)(\alpha_n^{-1}-z^2\psi^\prime_n(z))+\alpha_n^2z\psi^\prime_n(z)(\alpha_n^{-1}+z \psi_n(z))}{-\alpha_nz\psi_n(z)}\\
&=& \frac{\alpha_n^2\alpha_n^{-1}(\psi_n(z)+z\psi^\prime_n(z))}{-\alpha_nz\psi_n(z)}
 = - \frac{\psi_n(z)+z\psi^\prime_n(z)}{z\psi_n(z)}\\
& =-& \frac{1}{z}-\frac{\psi^\prime_n(z)}{\psi_n(z)} = -\frac{1}{z} - \frac{\sy^\prime\sy+2z\theta_n(z)+z^2\theta_n^\prime(z)}{1+z\sy^\prime\sy+z^2\theta_n(z)} \,.
\end{eqnarray*}
Finally, we derive the representation \eqref{thetrep} $\theta_n(z)$  using Woodbary matrix inversion lemma [see, e.g., \cite{hornjohn1985}]:
\begin{eqnarray*}
  \theta_n(z)&=&\sy^\prime(\bS_n^+-z\bI)^{-1}\sy\\
 &=& \sy^\prime\left(1/n\by_n(1/n\by^\prime_n\by_n)^{-2}\by^\prime_n-z\bI \right)^{-1}  \sy\\
&=& \sy^\prime \left(-\frac{1}{z}\bI-\frac{1}{z^2}1/n\by_n(\by^\prime_n\by_n)^{-1}(\bI-\frac{1}{z}(\by^\prime_n\by_n)^{-1})^{-1}(\by^\prime_n\by_n)^{-1}\by_n^\prime  \right)\sy\\
&=&-\frac{1}{z}\sy^\prime\sy-\frac{1}{z^2}\sy^\prime 1/n\by_n(1/n\by^\prime_n\by_n)^{-1}(\bI-\frac{1}{z}(1/n\by^\prime_n\by_n)^{-1})^{-1}(1/n\by^\prime_n\by_n)^{-1}\by_n^\prime  \sy\\
&=& -\frac{1}{z}\sy^\prime\sy-\frac{1}{z^2}\frac{1}{n}\bi_n^\prime (1/n\by_n^\prime\by_n)(1/n\by^\prime_n\by_n)^{-1}(\bI-\frac{1}{z}(1/n\by^\prime_n\by_n)^{-1})^{-1}(1/n\by^\prime_n\by_n)^{-1}(1/n\by_n^\prime\by_n)\bi_n \\
&=&  -\frac{1}{z}\sy^\prime\sy -\frac{1}{z^2}\frac{1}{n}\bi_n^\prime (\bI-\frac{1}{z}(1/n\by^\prime_n\by_n)^{-1})^{-1}\bi_n\\
&=&  -\frac{1}{z}\sy^\prime\sy +\frac{1}{z}\frac{1}{n}\bi_n^\prime ((1/n\by^\prime_n\by_n)^{-1}-z\bI)^{-1}\bi_n.
\end{eqnarray*}

\bigskip

\bigskip

{\bf Acknowledgements.} The authors would like to thank
M. Stein who typed this manuscript with considerable technical
expertise. The work of H. Dette was supported by the DFG Research Unit 1735, DE 502/26-2

\setlength{\bibsep}{4pt}
\bibliography{moorepenrose}

\end{document}